\documentclass[12pt,a4paper]{article}

\usepackage{graphicx}
\usepackage{amsmath}
\usepackage{latexsym}
\usepackage{amssymb}
\usepackage{hyperref}

\title{Catastrophe Management and Inter-Reserve Distance for Marine Reserve Networks}

\begin{document}

\author{
L.D. Wagner\footnote{Department of Mathematics, University of Queensland, QLD 4072, Australia and ARC Centre for Complex Systems, University of Queensland, QLD 4072, Australia. {\tt ldw@maths.uq.edu.au}}
\and 
J.V. Ross\footnote{Department of Mathematics, University of
Queensland, QLD 4072, Australia. {\tt jvr@maths.uq.edu.au}}
\and 
H.P. Possingham\footnote{School of Integrative Biology, University of Queensland, QLD 4072, Australia and Department of Mathematics, University of Queensland, QLD 4072, Australia. {\tt h.possingham@uq.edu.au}}
}

\maketitle

\abstract{We consider the optimal spacing between marine reserves for maximising the viability of a species occupying a reserve network. The closer the networks are placed together, the higher the probability of colonisation of an empty reserve by an occupied reserve, thus increasing population viability. However, the closer the networks are placed together, the higher the probability that a catastrophe will cause extinction of the species in both reserves, thus decreasing population viability. Using a simple discrete-time Markov chain model for the presence or absence of the species in each reserve we determine the distance between the two reserves which provides the optimal trade-off between these processes, resulting in maximum viability of the species.

}

\section{Introduction}

The design of reserve networks for marine conservation is a contentious issue 
from both theoretical and practical points of view.  One of the difficult theoretical issues in marine 
reserve system design concerns the optimal connectedness of the reserves within the network. The closer two reserves are placed together, the more likely it is that a population occupying one reserve will colonise the other, increasing population viability. However, most theories fail to consider the possibility of catastrophic events (Possingham et al. (2000), Shafer (2001), Sala et al. (2002) and Allison et al. (2003)); the closer two reserves are placed together the more likely they are to be struck by the same catastrophic event, thus decreasing viability. Consequently, there exists a natural tension between these two processes. Whilst this trade-off has been identified by many authors, a detailed theoretical investigation of its influence and the subsequent analytic determination of an optimal spacing, has, to the best of our knowledge, not been undertaken.

In this paper we explore the question of how closely two marine reserves should be placed in order to maximise the viability of a species in a reserve network under the competing processes of colonisation and catastrophe. We use a simple discrete-time Markov chain metapopulation model to explore the effect of reserve spacing on the viability of the metapopulation. The mathematical fomulation of the model follows closely the work of Day and Possingham (1995). We use the quasi-stationary distribution (the distribution of the metapopulation conditioned upon the population being extant) to determine an optimal inter-reserve spacing. The analysis relies on finding eigenvalues of the transition matrix of the Markov chain. The rate at which the metapopulation decays to extinction from quasi-equilibrium is given by the second eigenvalue of the probability transition matrix (Darroch and Senata (1965)). Thus, having found the second eigenvalue of the Markov chain, we can investigate the influence of inter-reserve distance on metapopulation viability, and by maximising this eigenvalue determine the optimal inter-reserve distance.

Having constructed a baseline model for the viability of a marine metapopulation, we will then examine how the inclusion of two further processes impacts optimal reserve spacing. The inclusion of an external recruitment process allows us to model a more expansive population which is not just confined to the reserves under consideration (Sweaerer et al. (2002) and Cowen et al. (2000)). The analysis of this new model is similar to that for the baseline model, however we must now use the stationary distribution to investigate the effect of external recruitment since extinction (absorption) is not possible. We determine the optimal spacing by minimising the probability of having both reserves empty in equilibrium. We extend the model further by including the possibility of local reserve extinction, in which the population can become extinct in the absence of catastrophes. These two additional processes allow us to create a more general and realistic model of a marine metapopulation and to develop a more comprehensive idea of optimal inter-reserve spacing. This is a significant extension over the recent work by McCarthy et al. (2005) and a related uncertainty analysis of that model (Halpern et al. (2006)).

In the baseline case, and the two extensions to this model, we find reserve spacings that maximise the viability of the metapopulation. Our analytical results are an advance on existing theory, which is largely driven by results from Monte Carlo simulation models and often ignores the effect of catastrophes.

\begin{figure}[h!] \label{diagram}
\begin{center}
\includegraphics[width=195pt]{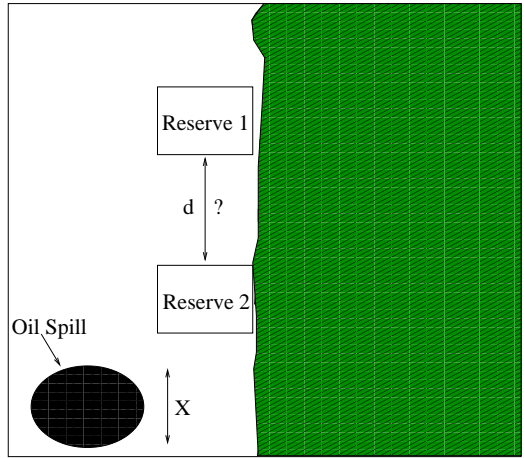}
\caption{Schematic of an example scenario where an oil spill has occurred along a coast line with two defined marine reserves which are placed at some distance apart.}
\end{center}
\end{figure}

\section{The Baseline Model}

Our baseline model consists of two processes: colonisation, and extinction caused by catastrophic events. This model is constructed to find the optimal spacing between two marine reserves under these competing processes. A catastrophe in this case is an extinction event which affects one or both populations. This occurs as either a direct reduction in population size or a process which degrades the habitat to such an extent that no organism will survive after a given time period, e.g. due to coral bleaching (Allison et al. (2003)).

For each reserve we model the presence or absence of the marine species within that reserve (Akcakaya and Ginzburg (1991)).  We also assume that both patches are identical, thus reducing the state space to three states since we only need to keep track of whether zero (total extinction), one or both patches are occupied. 

We have used basic probability theory to construct a discrete-time Markov chain for each of the processes -- colonisation and catastrophe -- considered in this model. To construct our probability transition matrix we need to determine the probabilities for each possible transition in our state space. Our model examines two marine reserves which lie at some point along a coastline. These two reserves are separated by a distance $d$ and are susceptible to being struck by some catastrophic event, such as an oil spill, which has a probability $r$ of occurrence in each time step. We assume that the probability that the catastrophe affects the second reserve decays exponentially with the distance $d$ between the reserves, and additionally depends upon the mean catastrophe size $\mu$. Based upon these assumptions we construct an extinction matrix $E$. An empty reserve may be colonised by an occupied reserve with a probability that decays exponentially with the distance $d$ between the reserves, and also depends upon a fitness parameter $\alpha$, which may be estimated for a particular marine system and accounts for the landscape between a particular reserve pair. The fitness parameter $\alpha$ may also be interpreted as the inverse of the average distance an individual disperses over the time period. Based upon these assumptions we construct a colonisation matrix $C$. Having constructed these matrices we can calculate the transition matrix $A$ for the complete process. We assume that extinction occurs before colonisation, and thus $A=E \times C$ (Gilpin (1992), Day and Possingham (1995)). This order of matrix multiplication is usually chosen since it is a ``worst-case" scenario, in that it will overestimate the probability of extinction. Additionally, in the present study we are only considering the second eigenvalue of the matrix $A$, and thus the order of multiplication does not affect the results. To calculate the probability of survival at some future time we may recursively evaluate $P_{t +1}=P_tA$, where $P_t$ is a vector containing the probabilities of being in each state at time step $t$.

\subsection{Colonisation Process}

We now need to establish the colonisation probability matrix. In this model we have that an unoccupied reserve can be colonised by an occupied reserve in each time step. Initially we construct a state space diagram (Figure (2)), which shows the four possible transitions. If both reserves are empty, then there are no individuals in the reserve network to colonise the reserves, and thus the zero state is clearly absorbing.

\begin{figure}[h!]
\begin{center}
\includegraphics[width=195pt]{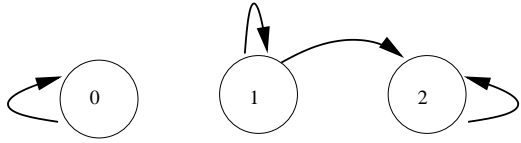}
\caption{State space with transitions between the numbers of occupied patches for the colonisation transition paths in our two patch metapopulation.}
\end{center}
\end{figure}

The process of migration from one patch to another is dependant upon the fitness $\alpha$ of a population within the landscape and the distance $d$ that the two reserves are placed apart from each other. This form of colonisation probability is ubiquitous in the literature and has been empirically observed for many metapopulations (Gilpin (1992), Hanski (1994a, 1994b) and Cowen et al. (2000)). Hence, the colonisation matrix is given by
\begin{equation}
C=
\begin{pmatrix}
	1 & 0 & 0\\
	0 & 1-e^{(-\alpha d)} & e^{(-\alpha d)}\\
	0 & 0 & 1
\end{pmatrix}.
\end{equation}

In Figure $(3)$ we can see that the chance one patch is colonised from another decays with distance and the relationship is dependent upon the fitness parameter $\alpha$. In a marine environment, the process of colonisation takes place continuously through time, implying that the dispersal of marine organisms in all stages of their life cycle, whether that is dispersal of larvae or the translocation of adults, plays an important role in metapopulation viability (Cowen et al. (2000)). Furthermore, organisms who are able to translocate as a stress response, related to some catastrophic event affecting their environment, will also migrate to an unaffected reserve in the manner proposed by $C$.

\begin{figure}[!h] \label{colon2}
\begin{center}
\scalebox{0.5}{
\includegraphics{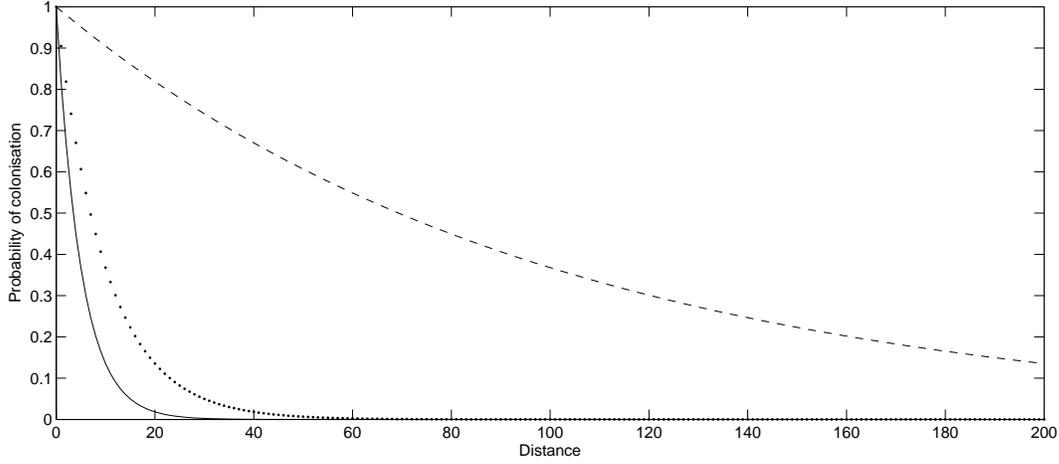}}
\caption{Graph of the probability of colonisation $\exp(-\alpha d)$ with respect to distance $d$, for $\alpha=0.01$ (solid line), $0.1$ (dots) and $0.2$ (dashed line).}
\end{center}
\end{figure}

\subsection{Extinction Process}

We now derive the extinction probability matrix. Once again we construct a state space diagram (Figure (4)), which shows the six possible transitions. Since we are only considering the process of extinction, an increase in the number of occupied reserves is impossible. 

A catastrophe that affects at least one reserve, arrives with probability $r$ within each time step. When a catastrophe occurs, we assume that it affects both reserves with probability $\exp(-d/\mu)$, and therefore affects only one reserve with probability $1-\exp(-d/\mu)$, where $\mu$ is the mean catastrophe size. This formulation assumes that the probability that a catastrophe affects both reserves decays exponentially with increasing inter-reserve distance. This is a reasonable assumption, in particular for catastrophes with a radius that is exponentially distributed in size, with an average radius of $\mu$.

\begin{figure}[h!]
\begin{center}
\includegraphics[width=190pt]{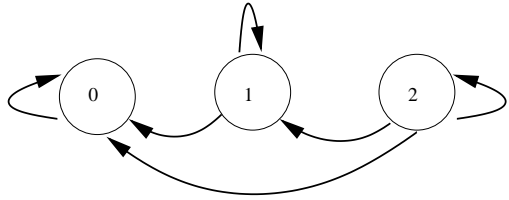}
\caption{State space with transitions between the number of occupied patches for the catastrophe transition paths in our two patch metapopulation.}
\end{center}
\label{fig:ss2}
\end{figure}

We shall now derive the functional forms for each possible transition. Firstly, if neither reserve is occupied, then the system will remain unoccupied since we are only considering catastrophic events. Thus, the entry $E(0,0)$ of our transition probability matrix $E$ will be $1$, and $E(0,1)$ and $E(0,2)$ will both be $0$. When both reserves are occupied, we have that the system will only remain fully occupied if a catastrophe does not occur within the time step, thus $E(2,2) = 1-r$. If a catastrophe does occur (with probability $r$), it removes only one reserve with probability $1-\exp(-d/\mu)$ and both reserves with probability $\exp(-d/\mu)$. Thus, $E(2,1)$ and $E(2,0)$ are $r(1-\exp(-d/\mu))$ and $r\exp(-d/\mu)$, respectively. Finally, when only one reserve is occupied, we have that the system will remain with only one reserve occupied if either a catastrophe does not occur (with probability $1-r$), or a catastrophe occurs but only affects the unoccupied reserve -- this occurs with probability $(1/2)(1-\exp(-d/\mu))$, since a catastrophe affects only one reserve with probability $(1-\exp(-d/\mu))$, and with probability $1/2$ this was the one occupied (since we are assuming that each patch is identical and that the catastrophe is equally likely to affect either patch initially). Thus, $E(1,1) = (1-r) + (r/2)(1-\exp(-d/\mu))$. The state of the system will change to no reserves occupied if the catastrophe affects both reserves, or affects only one reserve, with this reserve being the one that was initially occupied. Therefore, $E(1,0) = r[(1/2)(1-\exp(-d/\mu)) + \exp(-d/\mu)]$, and since we are only considering catastrophic events, $E(1,2) = 0$.

This generates our transition matrix for the catastrophe/extinction process,
\begin{equation} \label{transext}
E=
\begin{pmatrix}
	1 & 0 & 0\\
	\frac{r}{2}\left(1+e^{\frac{-d}{\mu}}\right) & 1 - \frac{r}{2}\left(1+e^{\frac{-d}{\mu}}\right) & 0\\
	re^{\frac{-d}{\mu}} & r\left(1-e^{\frac{-d}{\mu}}\right) & 1-r
\end{pmatrix}.
\end{equation} 

Figure~$(5)$ shows the probability $r\exp(-d/\mu)$ of a catastrophe occurring in a time step and it affecting both reserves, for different values of $\mu$ and $r=0.5$.

\begin{figure}[!h] \label{extint}
\begin{center}
\scalebox{0.5}{
\includegraphics{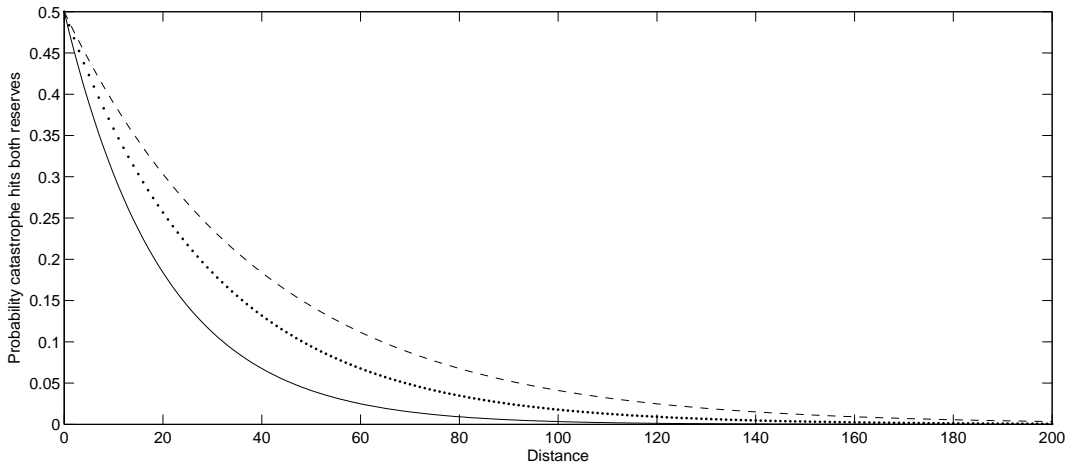}}
\caption{The probability of the extinction process hitting both reserves in a network within a year with respect to its mean size; $\mu=20$ (solid line), $30$ (dots) and $40$ (dashed line), all with respect to the occurrence of a catastrophe in one year $r=0.5$.}
\end{center}
\end{figure}

\section{Baseline Model Results}


As mentioned previously, to analyse the baseline model we use the second eigenvalue of the probability transition matrix. This eigenvalue provides the rate at which the metapopulation decays to extinction from quasi-stationarity. The graph below (Figure $(6)$) plots the second eigenvalue of $A = E\times C$ as a function of inter-reserve spacing $d$. The optimal inter-reserve distance is given by the value of $d$ for which the second eigenvalue is largest.


\begin{figure}[!h]
\begin{center}
\scalebox{0.5}{
\includegraphics{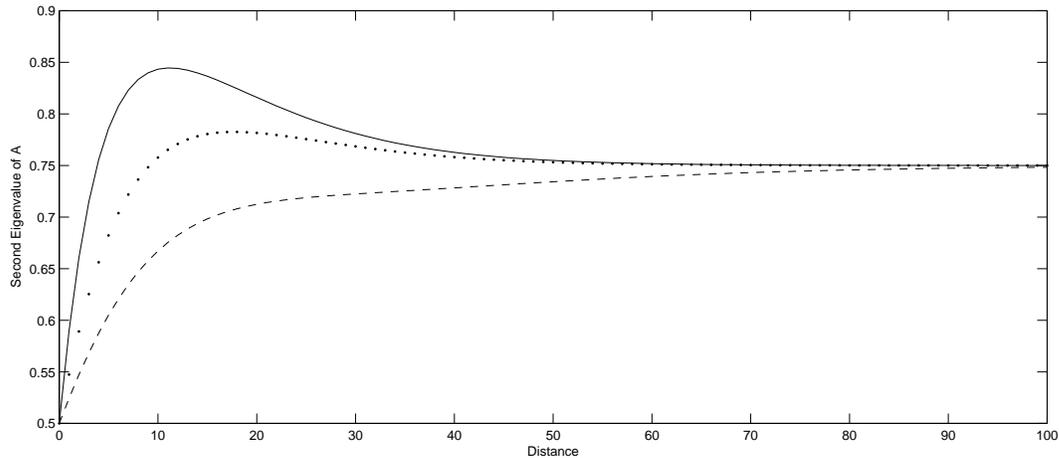}}
\caption{Second eigenvalue of $A$; $r=0.5$, $\alpha=0.1$ and $\mu=5$ (solid line), $10$ (dots), and $20$ (dashed line), as three different values for the mean size of disaster.}
\end{center}
\end{figure}

We can see from Figure $(6)$ that when the mean catastrophe size $\mu$ is smaller than, or comparable to, the mean dispersal distance $1/\alpha$, there exists a unique optimal inter-reserve spacing, signified by the hump in the curves corresponding to $\mu=5$ and $10$. It can also be seen that the optimal inter-reserve distance increases with increasing catastrophe size. As $\mu$ increases, eventually the trade-off between between colonisation and catastrophe disappears, and the optimal inter-reserve distance is found to be much larger. This emerges because as the distance increases (so as to attempt to avoid catastrophes affecting both reserves), the probability of colonisation becomes smaller. Eventually colonisation becomes essentially impossible, and as a consequence the optimal inter-reserve spacing is almost entirely determined by attempting to minimise the probability of a catastrophe affecting both reserves. In this situation, it is still optimal to place the reserves at the minimum distance for which the curve is no longer increasing (since catastrophes will usually be infrequent compared to colonisation events). This is to allow for the best chance of recolonisation.




Given particular parameter estimates, our model and method of analysis allows the calculation of the optimal inter-reserve distance. This generality is very helpful in attempting to determine the optimal placement of reserves within a complex marine park such as the Great Barrier Reef Marine Park, off Australia's east coast. Rather than having to use Monte Carlo simulations to estimate the probability of survival at a given distance in a two patch model such as this, we can calculate directly the effect of distance and plan more effectively for disaster.

\section{Model with External Recruitment}

In marine reserves, very few species have such a restricted occurrence that they only exist in two places. Thus, a logical extension of the model is to allow for external recruitment into the two patch metapopulation. This local patch effect, which has been coined the {\em rescue effect} by Brown and Kodric-Brown (1977), proposes that the migration of individuals into a struggling population may allow for the reduction in local extinction. The state space diagram in Figure (7), shows how the inclusion of external recruitment allows a locally extinct population (zero patches occupied), to now be recolonised (Jones et al. (1999) and Strathmann et al. (2002)). As such, the inclusion of external recruitment into our model creates a system without an absorbing state. However, we will refer to the state with both reserves unoccupied as being extinct.

\begin{figure}[ht!]
\begin{center}
\includegraphics[width=190pt]{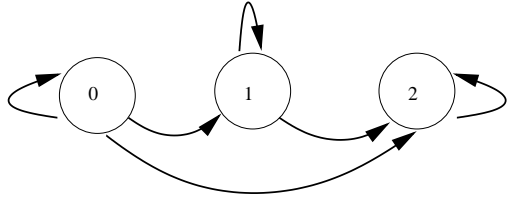}
\caption{State space with transitions between the number of occupied patches for the recruitment transitions paths in our two patch metapopulation.}
\end{center}
\end{figure}

To construct the new overall transition matrix $A$ we first need to construct an external recruitment matrix $R$, that tells us how external recruitment affects transitions between each state. The external recruitment matrix requires only one parameter $a$, being the probability that an empty patch is colonised by external recruits in a time step. Therefore we have that

\begin{equation}
R=
\begin{pmatrix}
	(1-a)^{2} & 2a(1-a) & a^{2}\\
	0 & 1-a & a\\
	0 & 0 & 1
\end{pmatrix}.
\end{equation}

The new transition matrix $A$, with catastrophes, colonisation and external recruitment combined is given by $A=E \times C \times R$. This order of matrix multiplication assumes that external recruitment occurs last, after the extinction and colonisation processes have occurred. The order of matrix multiplication does influence the results in this case, since we are interested in the first eigenvector (instead of the eigenvalue) of the matrix $A$. However, the order does not appear to have too much of an influence on the results from our investigations. This is partly due to the fact that we are only interested in finding the distance that minimises the probability of being extinct in equilibrium, rather than the actual probabilities of being in each of the states. The stationary distribution of $A$, which we obtain from the first eigenvector, tells us the probabilities of being in each of the states when the system is in equilibrium (Grimmett and Stirzaker (1992)). From the stationary distribution we can find an analytic solution for the placement of marine reserves in a network, by examining the first element within the stationary distribution; by taking one minus the first element within the stationary distribution we can calculate the probability that a species is extant in equilibrium. This is the probability which we need to maximize to effectively use this model to manage marine reserves. 

The graph of one minus the first element in the first eigenvector is presented below (Figure $(8)$). We can see that when the recruitment probability is small, we find similar results to those for the baseline model. The optimal inter-reserve distance does not appear to change too much following the introduction of external recruitment, with only a slight increase found. The optimal inter-reserve distance also appears to be relatively independent of the external recruitment probability. 

\begin{figure}[ht!]
\begin{center}
\scalebox{0.5}{
\includegraphics{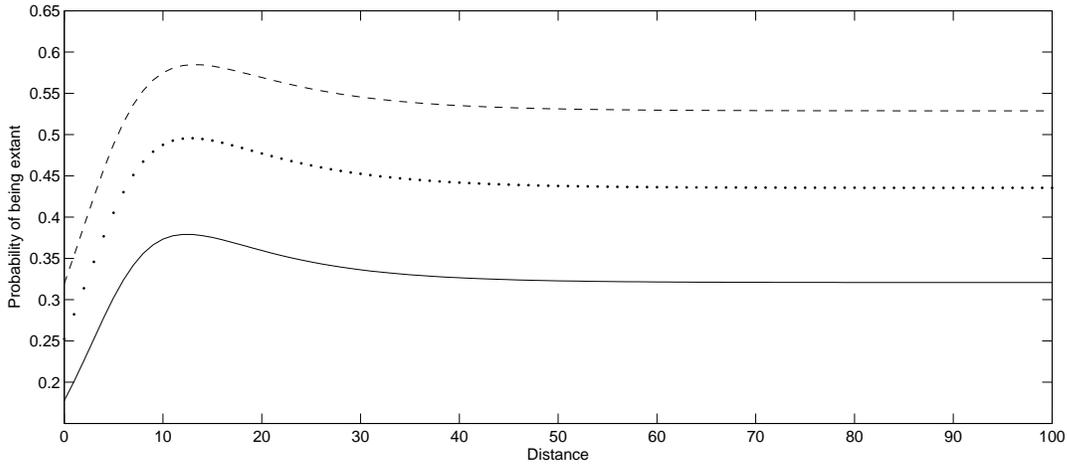}}
\caption{The probability of survival with external recruitment; $r=0.5$, $\alpha=0.1$, $\mu=5$ and $a=0.05$ (solid line), $0.075$ (dots) $0.10$ (dashed line), as three different local recruitment probabilities.}
\end{center}
\end{figure}

\section{Model with Local Extinction}

Another important aspect of metapopulations we consider in this investigation is the process of local extinction caused by non-catastrophic events (Akcakaya and Ginzburg (1991)). We construct a transition matrix $L$ which contains the probabilities of transitions between states due solely to local extinction events (e.g. which may occur through unfavourable local environmental conditions). The construction of $L$ is similar to that for the matrix $R$; the matrix $L$ only requires one parameter $b$, being the probability of local extinction within a time step. Therefore we have that

\begin{equation}
L=
\begin{pmatrix}
	1 & 0 & 0\\
	b & 1-b & 0\\
	b^{2} & 2b(1-b) & (1-b)^{2}
\end{pmatrix}
\end{equation}

The new transition matrix $A$, with catastrophes, colonisation, external recruitment and local extinction is given by $A=E \times L \times C \times R$. After constructing this Markov chain we perform the same analysis as with the previous model and find the spacing between reserves that maximizes the viability of the population (Figure $(9)$).

\begin{figure}[ht!]
\begin{center}
\scalebox{0.5}{
\includegraphics{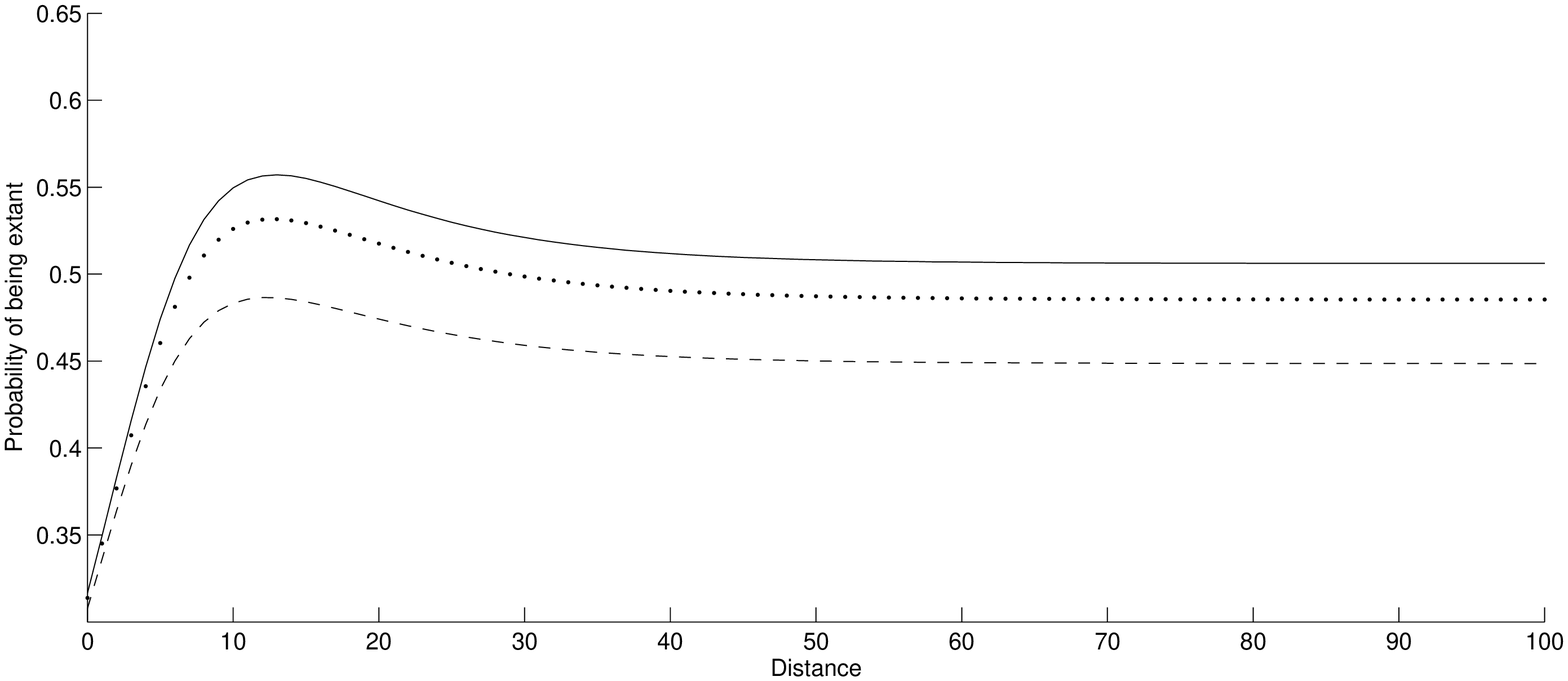}}
\caption{The probability of survival with local extinction; $r=0.5$, $\alpha=0.1$, $\mu=5$, $a=0.1$ and $b=0.025$ (solid line), $0.05$ (dots) and $0.1$ (dashed line), as three local extinction probabilities.}
\end{center}
\end{figure}

It can be seen that local patch extinction does not have too much of an effect on optimal inter-reserve spacing. The probability of being extant in equilibrium does however decrease, and thus the determination of the optimal inter-reserve distance becomes of more significance for ensuring maximum population persistence.

\newpage
\section{Discussion}

Ecologists and marine park managers have for a great deal of time relied upon practical experience and qualitative methods to design marine reserves for conservation purposes (Shafer (2001)). The additional inclusion of a variety of stake holders into the process of marine reserve design may also detract from the formulation of a conservation management plan that ensures the survival of a particular species. We have provided methods that allow the precise calculation of the optimal inter-reserve distance. This allows decision makers to have quick estimates for optimal inter-reserve spacing, and thus quantitative information for discussion with stake holders.


The results for our baseline model showed how to optimally space marine reserves under the competing processes of colonisation and catastrophe. This analysis involved maximising the second eigenvalue of the transition probability matrix, which results in the longest possible persistence of the species. We found that when the mean catastrophe size is smaller than, or comparable to, the mean dispersal distance, there exists a unique optimal inter-reserve spacing. As the mean catastrophe size increases, we arrive at a continuum of optimal inter-reserve distances. In this latter situation, the minimum of these optimal distances should usually be selected, to provide the maximum probability of colonisation.

We then included the possibility of external recruitment and local extinction events. With the expansion of the model there is now no absorbing state and the first element of the first eigenvector tells us the equilibrium probability that the species is absent. To calculate the optimal reserve design strategy we minimize this probability. We found similar behaviour to the situation where there is no external recruitment and local extinction, but the importance of determining the optimal inter-reserve spacing was highlighted.

As more information becomes available concerning species occupying marine landscapes, more accurate predictions for marine reserve design will be achievable. The values we have adopted in this paper appear typical for many species, with larval dispersal distances in the range of $10-100$ km~(Cowen et al. (2006)).

Our results are of most interest in situations were the process of dispersal and catastrophe are on a similar scale. The extension of our results to the multiple-reserve network will also be interesting, as catastrophe and colonisation will always be on a similar scale, at least for some subset of the whole network. In such situations it will also be of interest to look at combining our results with those of others concerning marine reserve design, see for example~McCarthy et al. (2005, 2006). In such situations there will always exist a trade-off between the size of the reserves, the number of reserves and the spacing, and the interaction between each with catastrophes will be of great importance (Airame et al. (2003)).

\bigskip
\noindent
{\textbf{Acknowledgements}}.
We would like to thank the referees for their valuable comments which did much to improve this paper. Additionally, LDW would like to acknowledge the support the Australian Research Council Centre for Complex Systems (ACCS) and JR wishes to acknowledge the support of the Australian Research Council Centre of Excellence for Mathematics and Statistics of Complex Systems (MASCOS). 


\begin{thebibliography}{199}

\bibitem{Airetal03} Airame, S., Dugan, J. E., Lafferty, K. D., Leslie, H., McArdle, D. A. and Warner, R. R., \textit{Applying ecological criteria to marine reserve design: A case study from the California Channel Islands.} Ecol. App., {\bf 13} (2003), pp.S170-S184

\bibitem{akc} Akcakaya, H.R. and Ferson, S., \textit{RAMAS/Space: Spatially Structured Population Models for Conservation Biology}, App. Biomath. (1992)

\bibitem{allison} Allison, G.W., Gaines, S.D., Lubchenco, J. and Possingham, H.P.,
\textit{Ensuring Persistence of Marine Reserves: Catastrophes Require Adopting
an Insurance Factor.} Ecol. App., {\bf 13} (2003), pp.S8-S24

\bibitem{brown} Brown, J.H. and Kodric-Brown, A., \textit{Turnover rates in insular biogeography: Effect of immigration on extinction.} Ecology, {\bf 58} (1977), pp.445-449

\bibitem{cowen} Cowen, R.K., Kamazima, M.M.L, Su, S., Paris, C.B. and Olsen, D.B, \textit{Connectivity of Marine Populations: Open or Closed?} Science, {\bf 287} (2000) pp.857-859   

\bibitem{cowen1} Cowen, R.K., Paris, C.B., Srinivasan, A., \textit{Scaling of Connectivity in Marine Populations}.
Science, {\bf 311} (2006), pp.522-527

\bibitem{Darroch65} Darroch, J. N. and Seneta, E. J. \textit{On quasi-stationary distributions in absorbing discrete-time finite Markov chains.} J. Appl. Probab. {\bf 2} (1965) pp.88-100 

\bibitem{Day95} Day, J.R. and Possingham, H.P.,  \textit{A stochastic metapopulation model with variability in patch size and position.} Theor. Pop. Biol. {\bf 48} (1995) pp.333-360

\bibitem{gilpin} Gilpin, M., \textit{Demographic stochasticity: A Markovian approach.}, J Theor. Biol. {\bf 154} (1992) pp.1-8

\bibitem{grim} Grimmett, G. and Stirzaker, D., \textit{Probability and Random Processes,} Oxford University Press, Oxford, 2nd Ed. (1992)

\bibitem{HRPM05} Halpern , B.S., Regan, H.M., Possingham, H.P. and McCarthy M.A. \textit{Accounting for uncertainty in marine reserve design.}, Ecology Letters {\bf 9} (2005) pp.2-11.

\bibitem{hanski} Hanski, I., \textit{A practical model of metapopulation population dynamics.} J. Animal Ecol. {\bf 63} (1994a) pp.151-162

\bibitem{hanskib} Hanski, I., \textit{Patch-occupancy dynamics in fragmented landscapes.} Trends Ecol. Evol. {\bf 9} (1994b) pp.131-135

\bibitem{jones} Jones, G.P., Milicich, M.J., Emlise, M.J. and Lunow, C., \textit{Self-Recruitment in a Coral Reef Population.}
Nature {\bf 402} (1999) pp.802-804

\bibitem{MTP05} McCarthy, M.A., Thompson, C.J. and Possingham, H.P., \textit{Theory for Designing Nature Reserves for Single Species.} Am. Nat. {\bf 165} (2005) pp.250-257

\bibitem{MTW06} McCarthy, M.A., Thompson, C.J. and Williams, N.S.G., \textit{Logic for designing natures reserves for multiple species.} Am. Nat. (To appear,(2006)).

\bibitem{possum1} Possingham, H.P., Ball., I.R. and Andelman, S., \textit{Mathematical methods for identifying representative reserve networks,} In: S. Ferson and M. Burgman (Eds) Quantitative methods for conservation biology. Springer-Verlag, New York, pp. 291-305 (2000)

\bibitem{Sala02} Sala, E., Aburto-Oropeza, O., Paredes, G., Parra, J., Barrera, J.C. and Dayton, P.K., \textit{A General Model  for Designing Networks of Marine Reserves.} Science {\bf 298} (2002) pp.1991-1993.

\bibitem{Shafer01} Shafer, C.L.,  \textit{Inter-reserve distance.} Biological Conservation {\bf 100} (2001) pp.215-227 


\bibitem{Stewart02} Stewart, R.R. and Possingham, H.P., \textit{A framework for systematic marine reserve design in South 
Australia: a case study.}, Inaugural World Congress on Aquatic Protected Areas, Cairns, August 2002

\bibitem{strath} Strathmann, R.R., Hughes, T.P., Kuris, A.M., Linedman, K.C., Morgan, S.G., Pandolfi, J.M. and Warner, R.R., \textit{Evolution of Local Recruitment and its Consquences for Marine Populations.} Bull. Marine Sci. {\bf 70} (2002), Suppl., pp.377-396

\bibitem{swear} Swearer, J.E., Shima, J.S., Hellberg, M.E., Thorrold, S.R., Jones, G.P., Robertson, D.R., Morgan, S.G., Selkoe, K.A., Ruiz, G.M. and Warner, R.R., \textit{Evidence of self recruitment in Demersal Marine Populations.} Bull. Marine Sci. {\bf 70} (2002) pp.251-271


\end{thebibliography}
\end{document}